\documentstyle[11pt,epic]{article}
\begin{document}
\newcommand{\up}{\vspace*{-0.2cm}}
\newcommand{\upp}{\vspace*{-0.3cm}}
\newcommand{\qed}{\hfill$\rule{.05in}{.1in}$\vspace{.3cm}}
\newcommand{\pf}{\noindent{\bf Proof: }}
\newtheorem{thm}{Theorem}
\newtheorem{lem}{Lemma}
\newtheorem{prop}{Proposition}
\newtheorem{prob}{Problem}
\newtheorem{ex}{Example}
\newtheorem{cor}{Corollary}
\newtheorem{conj}{Conjecture}
\newtheorem{cl}{Claim}
\newtheorem{df}{Definition}
\newtheorem{rem}{Remark}
\newcommand{\beq}{\begin{equation}}
\newcommand{\eeq}{\end{equation}}
\newcommand{\<}[1]{\left\langle{#1}\right\rangle}
\newcommand{\be}{\begin{enumerate}}
\newcommand{\ee}{\end{enumerate}}
\newcommand{\Bul}{\mbox{$\bullet$ } }
\newcommand{\al}{\alpha}
\newcommand{\ep}{\epsilon}
\newcommand{\si}{\sigma}
\newcommand{\om}{\omega}
\newcommand{\la}{\lambda}
\newcommand{\La}{\Lambda}
\newcommand{\Ga}{\Gamma}
\newcommand{\ga}{\gamma}
\newcommand{\im}{\Rightarrow}
\newcommand{\2}{\vspace{.2cm}}
\newcommand{\es}{\emptyset}

\title{\bf On Roman, Global and Restrained Domination in Graphs}

\author{A. Poghosyan and V. Zverovich\\
\\
{\footnotesize Department of Mathematics and Statistics}\\
{\footnotesize University of the West of England}\\
{\footnotesize Bristol, BS16 1QY}\\
{\footnotesize UK}\\
\\
}
\date{}

\maketitle

\begin{abstract}
In this paper, we present new upper bounds for the global
domination and Roman domination numbers and also prove that these
results are asymptotically best possible. Moreover, we give upper
bounds for the restrained domination and total restrained
domination numbers for large classes of graphs, and show that, for
almost all graphs, the restrained domination number is equal to
the domination number, and the total restrained domination number
is equal to the total domination number. A number of open problems
are posed.

\vspace*{.2cm} \noindent {\footnotesize {\it Keywords:} graphs,
Roman domination number, global domination number, restrained
domination number.}
\end{abstract}

\bigskip

\section{Introduction}
All graphs will be finite and undirected without loops and
multiple edges. If $G$ is a graph of order $n$, then
$V(G)=\{v_1,v_2,...,v_n\}$ is the set of vertices in $G$. Let
$N(x)$ denote the neighbourhood of a vertex $x$. Also let $
N(X)=\cup_{x\in X} {N(x)} $ and $ N[X]=N(X)\cup X.$ Denote by
$\delta(G)$ and $\Delta(G)$ the minimum and maximum degrees of
vertices of $G$, respectively. Put $\delta=\delta(G)$ and
$\Delta=\Delta(G)$.

A set $X$ is called a {\it dominating set} if every vertex not in
$X$ is adjacent to a vertex in $X$. The minimum cardinality of a
dominating set of $G$ is called the \emph{domination number}
$\gamma(G)$. The following fundamental result for the domination
number was proved by many authors \cite{alo,arn,lov1,pay}:

\begin{thm} [\cite{alo,arn,lov1,pay}]
\label{classic}
For any graph $G$,
$$
 \ga(G) \le {\ln(\delta+1)+1 \over \delta+1} n.
$$
\end{thm}

Let $H$ be a $k$-uniform hypergraph with $n$ vertices and $m$ edges.
The transversal number $\tau(H)$ of $H$ is the minimum cardinality of a set of vertices that
intersects all edges of $H$. Alon \cite{alo2} proved a fundamental result that if $k>1$, then
$$
\tau(H) \le {\ln k \over k} (n+m).
$$
He also showed that this bound is asymptotically best possible,
i.e. there exist $k$-uniform hypergraphs $H$ such that for
sufficiently large $k$,
$$
\tau(H) = {\ln k \over k} (n+m)(1+ o(1)).
$$
Alon \cite{alo2} gives an interesting probabilistic construction
of such a hypergraph $H$. In fact, $H$ is a random $k$-uniform
hypergraph on $[k\ln k]$ vertices with $k$ edges constructed by
choosing each edge randomly and independently according to a
uniform distribution on $k$-subsets of the vertex set. This
construction implies that the above bound for the domination
number is asymptotically best possible:

\begin{thm} [\cite{alo2}] \label{th2}
When $n$ is large there exists a graph $G$ such
that
$$
\ga(G) \ge {\ln(\delta+1)+1 \over \delta+1} n (1+o(1)).
$$
\end{thm}

 The concept of \emph{global domination} was
introduced by Brigham and Dutton \cite{dutton} and also by
Sampathkumar \cite{samp}. It is a variant of the domination
number. A set $X$ is called a \emph{global dominating set} if $X$
is a dominating set in both $G$ and its \emph{complement}
$\overline{G}$. The minimum cardinality of a global dominating set
of $G$ is called the \emph{global domination number}
$\gamma_{g}(G)$. There are a number of bounds on the global
domination number $\gamma_{g}(G)$. Brigham and Dutton in
\cite{dutton} give the following bounds on the global domination
number in terms of order, minimum and maximum degrees, and the
domination number of $G$:
\begin{thm}[\cite{dutton}]
If either $G$ or $\overline{G}$ is disconnected, then
$$
\ga_{g}(G)=\max \{ \gamma(G),\gamma(\overline{G})\}.
$$
\end{thm}

\begin{thm}[\cite{dutton}]
For any graph $G$, either $$ \ga_{g}(G)=\max \{
\gamma(G),\gamma(\overline{G})\} \,\, or \,\,\ga_{g}(G) \le \min
\{ \Delta(G),\Delta(\overline{G})\} + 1.$$
\end{thm}

\begin{thm}[\cite{dutton}]
For any graph $G$, if $\delta(G)= \delta(\overline{G}) \le 2$,
then
$$\ga_{g}(G) \le \delta(G)+2;$$ \\
 otherwise
 $$\ga_{g}(G) \le \max \{ \delta(G),\delta(\overline{G})\} + 1.$$
\end{thm}

 Another variant of the domination number, the \emph{Roman
domination number}, was introduced by Stewart \cite{ste}. In
\cite{ste} and \cite{ros}, a \emph{Roman dominating function
}(RDF) of a graph $G$ is defined as a function $f : V(G)
\rightarrow \{0, 1, 2\}$ satisfying the condition that every
vertex $u$ for which $f(u)= 0$ is adjacent to at least one vertex
$v$ for which $f(v) = 2$. The \emph{weight} of an RDF is defined
as the value $f(V(G)) = \sum_{v \in V(G)} f(v)$. The Roman
domination number of a graph $G$, denoted $\gamma_R(G)$, is equal
to the minimum weight of an RDF on $G$. In fact, Roman domination
is of both historical and mathematical interest. Emperor
Constantine had the requirement that an army or legion could be
sent from its home to defend a neighbouring location only if there
was a second army which would stay and protect the home. Thus,
there were two types of armies: stationary and travelling. Each
vertex with no army must have a neighbouring vertex with a
travelling army. Stationary armies then dominate their own
vertices, and a vertex with two armies is dominated by its
stationary army, and its open neighbourhood is dominated by the
travelling army. Thus, the definition of Roman domination has its
historical background and it can be used for the problems of this
type, which arise in military and commercial decision making. The
following results about the Roman domination number are known:

\begin{thm}[\cite{coc2}]
 For any graph $G$,
 $$
  \gamma(G) \le \gamma_R(G)\le 2\gamma(G).
 $$
 \end{thm}

\begin{thm}[\cite{coc2}]
For any graph $G$ of order $n$ and maximum degree $\Delta$,
 $$
 \gamma_R(G) \ge {2n \over {\Delta + 1}}.
 $$
\end{thm}

Telle and Proskurowski \cite{tel1} introduced restrained
domination as a vertex partitioning problem. A dominating set $X$
of a graph $G$ is called a {\em restrained dominating set} if
every vertex in $V(G)-X$ is adjacent to a vertex in $V(G)-X$. If,
in addition, every vertex of $X$ is adjacent to a vertex of $X$,
then $X$ is called a {\em total restrained dominating set}. The
minimum cardinality of a restrained dominating set of $G$ is the
\emph{restrained domination number} $\ga_r(G)$, and the minimum
cardinality of a total restrained dominating set of $G$ is the
\emph{total restrained domination number} $\ga_{\mathrm tr}(G)$.
For these parameters, the following upper bounds have been found:

\begin{thm} [\cite{dan1}]
If $\delta(G)\ge 2$, then
$$
\ga_r(G) \le n-\Delta.
$$
\end{thm}

\begin{thm} [\cite{hen1}]
If $G$ is a connected graph with $n\ge 4$, $\delta \ge 2$ and
$\Delta\le n-2$, then
$$
\ga_{\mathrm tr}(G) \le n- {\Delta \over 2} -1.
$$
\end{thm}

In this paper, we present new upper bounds for the global and
roman domination numbers, and show that our results are
asymptotically best possible. Moreover, we give upper bounds for
the restrained domination and total restrained domination numbers
for large classes of graphs. A number of open problems are posed.

\section{Upper Bounds for the Global Domination Number}

The following theorem provides an upper bound for the global
domination number. In what follows, we denote
$\bar{\delta}=\delta(\overline{G})$ and
$$
{\delta'} = \min \{\delta, \bar{\delta}\}.
$$

\begin{thm}\label{main1}
 For any graph $G$ with $\delta'>0$,
$$
\ga_{g}(G) \le \left(1-{{\delta'} \over {2}^{1/{\delta'}}
\,(1+{\delta'})^{1+1/{\delta'}}} \right) n.
$$
\end{thm}

\pf Let $A$ be a set formed by an independent choice of vertices
of $G$, where each vertex is selected with the probability
$$
p = 1-{1 \over {2}^{1/{\delta'}} \,(1+{\delta'})^{1/{\delta'}}}.
$$
Let us denote $ B = V(G)-N[A]$ and $ C = \{v_i \in V(G), v_i
\mbox{\,is not dominated by\,} A\,\mbox{in} \, \overline{G}\}.$ It
is easy to show that
\begin{eqnarray*}
{\mathbf P}[v_i\in B]&=&(1-p)^{1+\deg (v_i)} \, \,\,\,\,\,\,\,\, \,\,\,\,\,\, \le (1-p)^{1+\delta}\\
{\mathbf P}[v_i\in C]&=&(1-p)^{1+(n-\deg (v_i)-1)} \le
(1-p)^{1+\bar{\delta}}.
\end{eqnarray*}
 It is obvious that the set $ D =  A \cup B \cup C $ is a
global dominating set.
 The expectation
of $|D|$ is
\begin{eqnarray}
\nonumber
{\mathbf E}[|D|] &\le& {\mathbf E}[|A|]+{\mathbf E}[|B|]+{\mathbf E}[|C|]\\
\nonumber
&=& pn+ \sum_{i=1}^{n}{\mathbf P}[v_i\in B]+ \sum_{i=1}^{n}{\mathbf P}[v_i\in C]\\
\nonumber
&\le& pn + (1-p)^{1+\delta}n +(1-p)^{1+\bar{\delta}}n\\
\nonumber
&\le& pn + 2(1-p)^{1+\min \{\delta, \bar{\delta}\}}n\\
\nonumber &=& pn + 2(1-p)^{1+\delta'}n \\ \label{2} &=&
\left(1-{{\delta'} \over {2}^{1/{\delta'}}
\,(1+{\delta'})^{1+1/{\delta'}}} \right) n,
\end{eqnarray}
as required. The proof of the theorem is complete. \qed

The proof of Theorem \ref{main1} implies the following upper bound, which is asymptotically same as the bound of Theorem \ref{main1}.
\begin{cor}\label{Corollary 2}
For any graph $G$,
$$
\ga_{g}(G) \le {\ln({\delta'}+1) + \ln2+1 \over {\delta'}+1} n.
$$
\end{cor}

\pf Using the inequality $1-p \le e^{-p},$ we obtain the following
estimation of the expression (\ref{2}):
$$ {\mathbf E}[|D|] \le pn + 2e^{-p ({\delta'} +1)}
n. $$
 If we put $p=\min \lbrace 1, {\
 \ln({\delta'}+1)+\ln2 \over
{{\delta'}+1}}\rbrace$, then
$$ {\mathbf E}[|D|] \le {\ln({\delta'}+1) + \ln2+1 \over {\delta'}+1} n, $$ as required. \qed

We now prove that the upper bound of Corollary \ref{Corollary 2},
and therefore of Theorem \ref{main1}, is asymptotically best
possible.

\begin{thm}\label{best_possible}
When $n$ is large there exists a graph $G$ such that
$$
\ga_g(G) \ge {\ln(\delta'+1)+\ln2+1 \over \delta'+1} n (1+o(1)).
$$
\end{thm}

\pf Let us modify Alon's probabilistic construction described in
the introduction as follows. Let $F$ be a complete graph
$K_{[\delta\ln{\delta}]}$, and let us denote $F = V(F)$. Next, we
add a set of new vertices $V=\{v_1,...,v_{\delta}\}$, where each
vertex $v_i$ is adjacent to $\delta$ vertices that are randomly
chosen from the set $F$. Let us add a new component $K_{\delta+1}$
and denote the resulting graph by $G$, which has
$n=[\delta\ln{\delta}] + 2\delta+1$ vertices. Note that
$\delta'=\delta$ because $\bar{\delta}>\delta$. We will prove that
with high probability
$$
\ga_g(G) \ge  {\ln\delta' \over \delta'} n (1+o_{\delta'}(1)) =
{\ln\delta \over \delta} n (1+o_{\delta}(1)) =
\ln^{2}{\delta}(1+o_{\delta}(1)).
$$
Let us denote by $H$ the graph $G$ without the component
$K_{\delta+1}$. It is obvious that
$$\ga_g(G) =  \ga(H)+1.$$
Therefore, the result will follow if we can prove that  with high
probability
$$
\gamma(H)
> \ln^2\delta (1+o_{\delta}(1)).
$$

Without loss of generality we may only consider dominating sets in
$H$ that are subsets of $F$. Let us consider a dominating set $X$
in $H$ such that $X \subseteq F$ and $|X| \leq \ln^{2}{\delta}-
\ln\delta\ln{\ln^5 \delta}$. It is easy to show that the
probability of the set $X$ not dominating a vertex $v_i \in V$ is
\begin{eqnarray*}
{\mathbf P}[X \,{\rm{ does\, not\, dominate}}\,
v_i]&=&{{\pmatrix{|F| - |X| \cr \delta}} \over {\pmatrix{|F| \cr
\delta}}} \ge {\left(|F| - |X| - \delta \over |F| - \delta\right)}
^{\delta} = {\left(1-{|X| \over |F| - \delta}\right)} ^{\delta}.
\end{eqnarray*}
Using the inequality $1-x \ge e^{-x} (1 - x^2)$ if $x < 1$, we
obtain the following estimation:
\begin{eqnarray*}
{\mathbf P}[X \,{\rm{ does\, not\, dominate}}\, v_i]&\ge& e^{-
{\ln^2\delta - \ln\delta \ln{\ln^5\delta} \over \delta \ln\delta -
\delta}\delta}\left(1 - \left({\ln^2\delta - \ln\delta
\ln{\ln^5\delta} \over \delta \ln\delta -
\delta}\right)^2 \right)^{\delta}\\
&=& e^{-\ln\delta + \ln{\ln^5\delta} \over 1 - 1/\ln\delta}\left(1
+ o_{\delta}(1)\right)\\
&=& e^{\ln{({\ln^{5}\delta \over \delta})}(1 +
o_{\delta}(1))}\left(1
+ o_{\delta}(1)\right)\\
&=& \left({\ln^{5}\delta  \over \delta}\right)^{1+o_{\delta}(1)}(1
+
o_{\delta}(1))\\
&\ge& {\ln^4\delta \over \delta}.
\end{eqnarray*}
Thus, we conclude that
$$
{\mathbf P}\,[X \,{\rm dominates}\, V] \le \left(1- {\ln^4\delta
\over \delta}\right)^{\delta} \le e^{-\ln^4\delta}.
$$
It is obvious that the number of choices for the set $X$ is less
than $ \sum_{i=0}^{\ln^2\delta}{\pmatrix{|F| \cr i}}. $ We have
$${\sum_{i=0}^{\ln^2\delta}\pmatrix{|F| \cr
i}} <  \ln^2\delta \pmatrix{\delta\ln\delta \cr \ln^2\delta} <
(\delta\ln\delta)^{\ln^2\delta} < e^{2\ln^3\delta}.
$$
Now we can estimate the probability that the domination number of
the graph $H$ is less than or equal to $\ln^{2}{\delta}-
\ln\delta\ln{\ln^5 \delta}$:
$$
{\mathbf P}\left[\gamma(H) \le \ln^{2}{\delta}- \ln\delta\ln{\ln^5
\delta}\right] < \sum_{i=0}^{\ln^2\delta}{\pmatrix{|F| \cr
i}}{\mathbf P} \,[X \,{\rm dominates}\, V] < e^{2\ln^3\delta -
\ln^4\delta} = o_{\delta}(1).
$$
Therefore, with high probability $\gamma(H)
> \ln^{2}{\delta}-\ln\delta\ln{\ln^5 \delta} = \ln^2\delta (1+
o_{\delta}(1)),$ as required. The proof of the theorem is
complete. \qed

\section{Upper Bounds for the Roman Domination Number}

The following theorem provides an upper bound for the Roman
domination number:

\begin{thm}\label{main}
For any graph $G$ with $\delta>0$,
$$
\ga_{R}(G) \le 2\left(1-{{2^{1/{\delta}}\delta} \over
(1+{\delta})^{1+1/{\delta}}} \right) n.
$$
\end{thm}

\pf Let $A$ be a set formed by an independent choice of vertices
of $G$, where each vertex is selected with the probability
$$
p = 1-{\left({2} \over {1+{\delta}}\right)}^{1/{\delta}}.
$$
We denote $ B = N[A]-A$ and $ C = V(G)-N[A].$ Let us assume that
$f$ is a function $f : V(G) \rightarrow \{0, 1, 2\}$ and assign
$f(v_i)=2$ for each $v_i \in A$, $f(v_i)=0$ for each $v_i \in B$
and $f(v_i)=1$ for each $v_i \in C$. It is obvious that $f$ is a
Roman dominating function and $f(V(G))=2|A|+|C|$.

It is easy to show that
\begin{eqnarray*}
{\mathbf P}[v_i\in C]&=&(1-p)^{1+\deg(v_i)} \,  \le
(1-p)^{1+\delta}.
\end{eqnarray*}
 The expectation
of $f(V(G))$ is
\begin{eqnarray}
\nonumber
{\mathbf E}[f(V(G))] &\le& 2{\mathbf E}[|A|]+ {\mathbf E}[|C|]\\
\nonumber
&=& 2pn+ \sum_{i=1}^{n}{\mathbf P}[v_i\in C]\\
\nonumber &\le& 2pn + (1-p)^{1+\delta}n\\ \label{roman} &=&
2\left(1-{\delta \,{ 2^{1/{\delta}}} \over
(1+{\delta})^{1+1/{\delta}}} \right) n.
\end{eqnarray}
Since the expectation is an average value, there exists a
particular Roman dominating function of the above order, as
required. The proof of the theorem is complete. \qed
\newline
Theorem \ref{main} implies the following upper bound.

\begin{cor}\label{Corollary 1}
For any graph $G$ with $\delta>0$,
$$
\ga_{R}(G) \le {2\ln({\delta}+1) - \ln4+2 \over {\delta}+1} n.
$$
\end{cor}

\pf Using the inequality $1-p \le e^{-p},$ we obtain the following
estimation of the expression (\ref{roman}):
$$ {\mathbf E}[f(V(G))] \le 2pn + e^{-p ({\delta} +1)}
n. $$
 If we put $p= {
 \ln({\delta}+1)-\ln2 \over
{{\delta}+1}}$, then
$$ {\mathbf E}[f(V(G))] \le {2\ln({\delta}+1) - \ln4 +2 \over {\delta}+1} n, $$ as required. \qed
\newline
Note that the result of Corollary \ref{Corollary 1} was also
proved in \cite{coc2}, even though the upper bound in \cite{coc2}
contains a misprint.

Now let us prove that the upper bound of Corollary \ref{Corollary
1}, and therefore of Theorem \ref{main}, is asymptotically best
possible.

\begin{thm}\label{best_possible1}
When $n$ is large there exists a graph $G$ such that
$$
\ga_R(G) \ge {2\ln(\delta+1)-\ln4+2 \over \delta+1} n (1+o(1)).
$$
\end{thm}

\pf Let $F$ be a complete graph $K_{[\delta\ln{\delta}]}$, and let
us denote $F = V(F)$. Next, we add a set of new vertices
$V=\{v_1,...,v_{\delta}\}$, where each vertex $v_i$ is adjacent to
$\delta$ vertices that are randomly chosen from the set $F$. The
resulting graph is denoted by $G$ and it has $n =
[\delta\ln{\delta}] + \delta$ vertices. We will prove that with positive probability
$$ \ga_R(G) \ge {2\ln\delta \over \delta} n (1+o_{\delta}(1))
= 2\ln^{2}{\delta}(1+o_{\delta}(1)).
$$Let $f = (D_0, D_1, D_2)$
be a $\ga_R$-function of $G$, i.e.$f$ is a Roman dominating
function and $f(V(G)) = \ga_R(G)$. It is easy to see that we may
assume that $D_2 \subseteq F$ and $D_1 \subseteq V$.\\

Let us consider two cases. If $|D_2|> \ln^2\delta -
\ln\delta\ln\ln^4\delta,$ then $f(V(G))
> 2\ln^2\delta(1+o_{\delta}(1))$, as required.
If $|D_2| \le \ln^2\delta - \ln\delta\ln\ln^4\delta,$ then, similar to the proof of Theorem \ref{best_possible}, we can
prove that the probability of the set $D_2$ dominating a vertex
$v_i \in V$ is
\begin{eqnarray}\label{exp_proof}
\nonumber {\mathbf P}[D_2 \,{\rm {dominates}}\, v_i] \le 1 -
{\ln^3\delta\over\delta}. \nonumber
\end{eqnarray}
Let us consider the random variable $|N(D_2) \cap V|$. The expectation of $|N(D_2) \cap V|$ is
$${\mathbf E}[|N(D_2) \cap V|]=\sum_{i=1}^{\delta}{\mathbf P}[D_2 \,{\rm
{dominates}}\, v_i] \le \delta - \ln^3\delta.$$ Thus we can
conclude that there exists a graph $G$, for which $|D_1| \ge
\ln^3\delta,$ i.e. $f(V(G)) \ge \ln^3\delta >
2\ln^{2}{\delta}(1+o_{\delta}(1))$, as required.  \qed


\section{Restrained and Total Restrained Domination}

Theorem \ref{classic} implies that when $\delta(G)$ is large,
$\ga(G)/n$ is close to 0 for any graph $G$. Similar results were
proved for the global and Roman domination numbers in the previous
sections. However, for the total restrained domination numbers
this is not the case, because for any $\delta$ there exists (see
\cite{hen1}) an infinite family of graphs $G$ with minimum degree
$\delta$, for which $\ga_{\mathrm tr}(G)/ n \rightarrow 1 $
 when $n$ tends to $\infty$. The above is also true for the
restrained domination number. Thus, for the class of all graphs, it
is impossible to provide an upper bound for these parameters similar
to the result of Theorem \ref{classic}. In this section, we will
give such upper bounds for large classes of graphs.

Let us first find the restrained domination number of a `typical'
graph. Let $0<p<1$ be fixed and put $q=1-p$. Denote by ${\cal
G}(n,{\bf P}[edge]=p)$ the discrete probability space consisting
of all graphs with $n$ fixed and labelled vertices, in which the
probability of each graph with $M$ edges is $p^{M}q^{N-M}$, where
$N= \pmatrix{n\cr2}$. Equivalently,  the  edges of a labelled
random graph are chosen independently and with the same
probability $p$. We say that a random  graph ${\bf G}$ {\it
satisfies a property} $Q$ if
$$
{\mathbf P}[{\bf G} \hbox{ has } Q] \to 1 \hbox{ as } n\to \infty
.
$$
If a random  graph ${\bf G}$ has a property $Q$, then we also say
that almost all graphs satisfy $Q$.

It turns out that, for almost all graphs, the restrained
domination number is equal to the domination number, which has two
points of concentration, and the total restrained domination
number is equal to the total domination number. This is formulated
in the following theorem, which is based on the fundamental
results of Bollob\'{a}s \cite{bol1} and Weber \cite{web1}. Remind
that a dominating set $X$ is called a {\em total dominating set}
if every vertex of $X$ is adjacent to a vertex of $X$. The {total
domination number} $\ga_t(G)$, which is one of the basic
domination parameters, is the minimum cardinality of a total
dominating set of $G$.

\begin{thm}
\label{almost}
For almost every graph, $\ga_r(G) = \ga(G)$ and
$\ga_{\mathrm tr}(G) = \ga_t(G)$. Moreover,
$$
\ga_r(G) = \lfloor\log n - 2\log\log n + \log\log e\rfloor +
\epsilon,
$$
where $\epsilon=1$ or $2$, and $\log$ denotes the logarithm with
base $1/q$.
\end{thm}

\pf Bollob\'{a}s \cite{bol1} proved that a random graph ${\bf G}$
satisfies
$$
\mid \delta ({\bf G}) - pn + (2pqn \log n)^{1/2} - {\Big( \frac
{pqn}{8\log n} \Big) }^{1/2} \log \log n\mid \; \le \, C(n){\Big(
\frac {n}{\log n} \Big) }^{1/2},
$$
where $C(n)\to \infty $ arbitrarily slowly. Therefore,
$$
\delta ({\bf
G}) = pn (1+o(1)).
$$

Weber \cite{web1} showed that the domination number of a random
graph ${\bf G}$ is equal to
$$
k+1 \quad \mbox{or} \quad k+2,
$$
where
$$
k=\lfloor\log n - 2\log\log n + \log\log e\rfloor
$$
and $\log$ denotes the logarithm with base $1/q$. Let us consider
a minimum dominating set $D$ of this size. We have
$$
|D|= \log n (1+o(1)).
$$
For any vertex $v\in V({\bf G})-D$ and large $n$,
$$
\deg v \ge \delta = pn(1+o(1))> \log n (1+o(1)) = |D|,
$$
since $p$ is fixed. Therefore, the vertex $v$ is adjacent to a
vertex in $V({\bf G})-D$, i.e. $D$ is a restrained dominating set.

Now let us consider a minimum total dominating set $T$, i.e.
$|T|=\ga_t({\bf G})$. It is not difficult to see that
$$
\ga_t({\bf G}) \le 2\ga({\bf G}).
$$
Therefore,
$$
|T| \le 2|D| = 2\log n (1+o(1)).
$$
Thus, for any vertex $v\in V({\bf G})-T$ and large $n$,
$$
\deg v \ge \delta = pn(1+o(1))> 2\log n (1+o(1)) \ge |T|,
$$
since $p$ is fixed. Therefore, the vertex $v$ is adjacent to a
vertex in $V({\bf G})-T$, i.e. $T$ is a total restrained
dominating set, which is also minimum. The result follows. \qed

However, the property of a `typical' graph stated in the above
theorem cannot be used as a bound for the (total) restrained
domination number for a given graph.  Let us find such upper
bounds for large classes of graphs.

\begin{prop} \label{pr2}
If $\delta>0$ and $n<\delta^2/(\ln\delta+1)$, then
$$
 \ga_r(G) \le {\ln(\delta+1)+1 \over \delta+1} n
$$
and
$$
 \ga_{\mathrm tr}(G) \le {\ln\delta +1 \over \delta} n.
$$
\end{prop}

\pf Using Theorem \ref{classic}, let us consider a dominating set
$D$ such that
$$
 |D| \le {\ln(\delta+1)+1 \over \delta+1} n.
$$
Note that the condition $n<\delta^2/(\ln\delta+1)$ can be written as
follows:
$$
\delta > {\ln\delta +1 \over \delta}n.
$$
Now, for any vertex $v\in V(G)-D$,
$$
\deg v \ge \delta > {\ln\delta +1 \over \delta}n >{\ln(\delta+1)+1
\over \delta+1} n \ge |D|.
$$
Therefore, the vertex $v$ is adjacent to a vertex in $V(G)-D$,
i.e. $D$ is a restrained dominating set.

Using the probabilistic method of the proof of Theorem
\ref{classic}, we can show that for any graph $G$ with $\delta>0$,
$$
\ga_t(G) \le {\ln\delta +1 \over \delta} n.
$$
Let us consider a total dominating set $T$ such that
$$
|T| \le {\ln\delta +1 \over \delta} n.
$$
For any vertex $v\in V(G)-T$,
$$
\deg v \ge \delta > {\ln\delta +1 \over \delta}n \ge |T|.
$$
Therefore, the vertex $v$ is adjacent to a vertex in $V(G)-T$,
i.e. $T$ is a total restrained dominating set. \qed

Note that the result of Bollob\'{a}s \cite{bol1} on the minimum
degree implies that the condition $n<\delta^2/(\ln\delta+1)$ is
satisfied for almost all graphs, i.e. Proposition \ref{pr2} gives
upper bounds for a very large class of graphs. Moreover, in the
class of graphs with $n<\delta^2/(\ln\delta+1)$, the upper bounds
of Proposition \ref{pr2} cannot be improved. This can be proved in
the same way as Theorem \ref{th2}.

The {\em matching number} of a graph $G$, denoted by $\beta_1(G)$,
is the largest number of pairwise non-adjacent edges in $G$. This
number is also called the {\em edge independence number}.

\begin{thm} \label{beta}
For any graph $G$ with $\delta>0$,
$$
 \ga_r(G) \le {2\ln(\delta+1)+\delta + 3 \over \delta+1}n -
 2\beta_1
$$
and
$$
 \ga_{\mathrm tr}(G) \le {2\ln\delta +\delta + 2 \over \delta}n -
 2\beta_1.
$$
\end{thm}

\pf Let us consider a minimum dominating set $|D|$ of the graph
$G$, i.e. $|D|=\ga(G)$. Let $M$ be a matching with $\beta_1(G)$
edges:
$$
M=(e_1,e_2,...,e_{\beta_1}).
$$
Without loss of generality we may assume that the first $k$ edges
of $M$ have at least one end in $D$, thus $\beta_1-k$ edges of $M$
have both ends in $V(G)-D$. It is obvious that
$$
k \le |D| = \ga(G).
$$
Therefore, at least $\beta_1(G) - \ga(G)$ edges in $M$ have both
end vertices in $V(G)-D$. Note that $\beta_1(G) \ge \ga(G)$,
because each vertex of a total dominating set $S$ has a private
neighbour not in $S$, thus providing a matching of size $|S|$,
which is at least $\ga(G)$.

Now we form a restrained dominating set $D'$ by adding to $D$ all
vertices not belonging to the last $\beta_1-k$ edges of $M$. We
obtain
$$
 \ga_r(G) \le |D'| = n- 2(\beta_1-k)\le n - 2\beta_1+2\ga.
$$
By Theorem \ref{classic},
$$
\ga(G) \le {\ln(\delta+1)+1 \over \delta+1} n.
$$
Therefore,
$$
 \ga_r(G) \le {2\ln(\delta+1)+\delta + 3 \over \delta+1}n -
 2\beta_1,
$$
as required.

Let us prove the latter upper bound. Consider a minimum total
dominating set $T$ and the above matching $M$. Using a similar
technique, we can construct a total restrained dominating set $T'$
such that
$$
 \ga_{\mathrm tr}(G) \le |T'| \le n - 2\beta_1+2\ga_t.
$$
Using the probabilistic method of the proof of Theorem
\ref{classic}, we can show that for any graph $G$ with $\delta>0$,
$$
\ga_t(G) \le {\ln\delta +1 \over \delta} n.
$$
Therefore,
$$
 \ga_{\mathrm tr}(G) \le {2\ln\delta +\delta + 2 \over \delta}n -
 2\beta_1,
$$
as required. \qed

A matching is called {\em perfect} if it contains all vertices of
a graph (or all vertices but one if $n$ is odd). The following
corollary follows immediately from the above theorem:

\begin{cor} \label{perfect matching}
If $G$ has a perfect matching, then
$$
 \ga_r(G) \le {\ln(\delta+1)+1 \over \delta+1} 2n + \epsilon
$$
and
$$
 \ga_{\mathrm tr}(G) \le {\ln\delta +1 \over \delta} 2n + \epsilon,
$$
where $\epsilon=0$ if $n$ is even and $\epsilon=1$ otherwise.
\end{cor}

It may be pointed out that the class of graphs with a perfect
matching includes all Hamiltonian graphs. It is well known that
almost all graphs are Hamiltonian \cite{moo1}, thus Corollary
\ref{perfect matching} provides upper bounds for a very large
class of graphs.

\section{Concluding Remarks and Open Problems}

By Theorem \ref{almost}, the total restrained domination number is
equal to the total domination number for almost every graph.
However, we do not know exact values of the total domination
number for almost all graphs. Such a result for the domination
number is known \cite{web1}.

\begin{prob}
For almost all graphs, find points of concentration of the total,
global and Roman domination numbers.
\end{prob}

Theorem \ref{classic} is formulated for all graphs and it gives an
excellent upper bound if $\delta$ is big. However,  for small
values of $\delta$, better (sharp) bounds are known:

\begin{thm} [Ore]
If $\delta(G) \ge 1$, then
$$
\ga(G)\le {n\over 2}.
$$
\end{thm}

\begin{thm} [\cite{mcc1}]
If $G$ is a connected graph with $\delta \ge 2$ and it is not
isomorphic to one of seven graphs (not shown here), then
$$
\ga(G)\le {2\over 5} n.
$$
\end{thm}

\begin{thm} [\cite{ree1}]
If $G$ is a connected graph with $\delta \ge 3$, then
$$
\ga(G)\le {3\over 8} n.
$$
\end{thm}

The above situation is also true for many upper bounds proved in
this paper. They are good when $\delta$ is not small. Can better
upper bounds be found for small values of $\delta$?

\begin{prob}
Determine sharp upper bounds for the global and Roman domination
numbers of a graph with small minimum degree.
\end{prob}

\begin{prob}
Determine sharp upper bounds for the restrained and total
restrained domination numbers of a graph with a perfect matching
and small minimum degree.
\end{prob}


\end{document}